\documentclass[11pt,a4paper,twoside]{amsart}
\usepackage{amsmath}
\usepackage{amssymb}
\usepackage{latexsym}

\newcommand{\co}{\colon\thinspace}    
\newcommand{\fnote}[1]{\footnote{\small sharp1}}
\newcommand{\inv}{^{-1}}              

\newcommand{\N}{{\mathbb N}}
\newcommand{\Z}{{\mathbb Z}}
\newcommand{\R}{{\mathbb R}}
\newcommand{\Q}{{\mathbb Q}}

\newcommand{\T}{{\mathbb T}}

\newcommand{\spt}{\mbox{supp}}
\newcommand{\Azero}{\mathcal{A}_0}

\newtheorem{theorem}{Theorem}
\newtheorem{proposition}[theorem]{Proposition}
\newtheorem{corollary}[theorem]{Corollary}
\newtheorem{definition}[theorem]{Definition}
\newtheorem{lemma}[theorem]{Lemma}
\newtheorem{remark}[theorem]{Remark}

\title{Vertices of Mather's beta function, II}
\author{Daniel Massart}
\date{\today}
\begin{document}

\begin{abstract}
If the $\beta$-function of a time-periodic Lagrangian on a manifold $M$ has a vertex at a $k$-irrational homology class $h$, then  $2k \leq \dim M$. Furthermore if $\dim M =2$ $h$ is rational.
\end{abstract}
\maketitle

\section{Introduction}
This paper addresses the problem of the differentiability of Mather's $\beta$-function for time-periodic Lagrangian systems. It started as an attempt to understand the result of \cite{Osvaldo} and reformulate it in the language of \cite{ijm}. Its aim is to unify under one method the various existing  results (\cite{Mather90}, \cite{Bangert94}, \cite{BIK}, \cite{ijm}, \cite{Osvaldo}). It fails in this respect since it does not contain Corollary 3 of \cite{ijm}, but it does improve a little bit on Theorem 1 of \cite{BIK} and \cite{Osvaldo}. 

The setting is as in \cite{Mather91} : $M$ is a closed, connected manifold. A time-periodic Lagrangian on $M$ is a $C^{2}$ function on $TM\times \T^{1}$, $\T^{1}$ being the unit circle, such that $L$ is convex and superlinear when restricted to the fibers of $TM$. An example to keep in mind is the sum of a Riemann  metric, viewed as a quadratic function on $TM$, and a time-periodic potential (a function on $M\times \T^{1}$).  The Euler-Lagrange equation gives rise to a flow $\Phi_t$ on $TM\times \T^{1}$. We make the additional assumption that  $\Phi_t$ is complete. See \cite{Fathi}, \cite{Bernard_pg}  for more background and references.

 Define $\mathcal{M}_{inv}$
to be the set of $\Phi_t$-invariant, compactly supported, Borel probability measures on $TM\times \T^{1}$.
Mather showed that the function (called action of the Lagrangian on measures)
	\[
	\begin{array}{rcl}
\mathcal{M}_{inv} & \longrightarrow & \R \\
\mu & \longmapsto & \int_{TM\times \T^{1}}	Ld\mu
\end{array}
\]
is well defined and has a minimum.  A measure achieving the minimum is called $L$-minimizing.

When $M=\T^{1}$, by Mather's Graph Theorem (\cite{Mather91}) an invariant measure can be given a rotation number just like an invariant measure of a circle homeomorphism. For other manifolds Mather proposed in \cite{Mather91} the following generalization. First he observed that if $\omega$ is a closed one-form on $M$ and $\mu \in \mathcal{M}_{inv}$ then the integral $\int_{TM\times \T^{1}}	\omega d\mu$ is well defined, and only depends on the cohomology class of $\omega$. By duality this endows  $\mu$ with a homology class : $\left[\mu\right]$ is the unique $h \in H_1 (M,\R)$ such that  
	\[
<h,\left[\omega \right]>= \int_{TM\times \T^{1}}	\omega d\mu 
\]
for any closed one-form $\omega$ on $M$. Besides, for any $h \in H_1 (M,\R)$, the set 
	\[ \mathcal{M}_{h,inv}:= \left\{\mu \in \mathcal{M}_{inv} \co \left[\mu\right]=h\right\}
\]
is not empty. Again the action of the Lagrangian on this smaller set of measures has a minimum, which is a function of $h$, called the $\beta$-function of the system. A measure achieving the minimum is called $(L,h)$-minimizing. 

There is a dual construction : if $\omega$ is a closed one-form on  $M$, then $L-\omega$ is a Lagrangian to which Mather's theory applies, and which has the same Euler-lagrange flow as $L$. Then the minimum over $\mathcal{M}_{inv}$ of $\int (L-\omega)d\mu$ is actually a function of the cohomology class of $\omega$, the opposite of which is called the $\alpha$-function of the system. An  $(L-\omega)$-minimizing measure is also called $(L,\omega)$-minimizing or $(L,c)$-minimizing if $c$ is the cohomology of $\omega$. In formal terms we have defined
	\[
	\begin{array}{rcl}
\beta \co H_1 (M,\R) & \longrightarrow & \R \\
h & \longmapsto & 
\min \left\{\int_{TM\times \T^{1}}	Ld\mu \co \left[\mu\right]=h\right\}\\
\alpha \co H^1 (M,\R) & \longrightarrow & \R \\
c & \longmapsto & 
\min \left\{\int_{TM\times \T^{1}}	(L-\omega)d\mu \co \left[\omega\right]=c\right\}.
\end{array}
\]

Mather proved that $\alpha$ and $\beta$ are convex, superlinear, and Fenchel dual of one another. The main geometric features of a convex function are its smoothness and strict convexity, or lack thereof. In the present setting they turn out to have interesting dynamical meanings as well. The prototype of all theorems in the subject is 
\begin{theorem}[\cite{Mather90}]\label{Mather}
If $M=\T^{1}$ then $\beta$ is differentiable at every irrational homology class. It is differentiable at a rational homology class if and only if periodic orbits in this class fill up $\T^{1}$. 
\end{theorem}
Since $H_1 (\T^{1},\R)=\R$ the word rational is self-explanatory. For other manifolds we need a bit of terminology. The torsion-free part of $H_1 (M,\Z)$ embeds as a lattice $\Gamma$ in $H_1 (M,\R)$. A class 
$h \in H_1 (M,\R)$ is called integer if it lies in $\Gamma$, and rational if $nh \in \Gamma$ for some $n \in \Z$. A subspace of $H_1 (M,\R)$ is called integer if it is generated by integer classes.

A convex function has a tangent cone at every point. We say it has a vertex at $x$ if its tangent cone at $x$ contains no straight line. An enticing question, suggested by Mather's theorem, is whether or not vertices of $\beta$  only occur at rational homology classes. As yet the best known result is
\begin{theorem}[\cite{Osvaldo}]
If the $\beta$-function of a Lagrangian $L$ has a vertex at a homology class $h$, then either $h$ is rational, or the support of every $(L,h)$-minimizing measure has Hausdorff dimension $\geq 3$.
\end{theorem}
Note that this contains  Theorem \ref{Mather} since by Mather's Graph Theorem (\cite{Mather91}) the support of a minimizing measure may be viewed as a subset of $M\times \T^{1}$ ; and if $ M= \T^{1}$,  $M\times \T^{1}$ cannot contain a subset of dimension three.  

We need to give a quantitative meaning to the irrationality of a homology class. The quotient $H_1 (M,\R)/\Gamma$ is a torus $\T^{b}$, where $b$ is the first Betti number of $M$. For $h$ in $H_1 (M,\R)$, the image of $\Z h$ in $\T^{b}$ is a subgroup of $\T^{b}$, hence its closure is a finite union of tori of equal dimension. This dimension is called the irrationality of $h$. It is zero if $h$ is rational. We say a class $h$ is completely irrational if its irrationality is maximal, i.e. equals $b$. Similar definitions can be made for vectors in $\R^{n}$ with the integer lattice $\Z^{n}$. Note that the irrationality of $h$ equals that of $Nh$ for $N\in \Z$, $N \neq 0$.
Our main contribution is
\begin{theorem}\label{vertices}
If the $\beta$-function of a Lagrangian $L$ on a manifold $M$ has a vertex at a homology class $h$, and $p$ is the irrationality of $h$, then the support of every $(L,h)$-minimizing measure has Hausdorff dimension $\geq 2p+1$. In particular $2p \leq \dim M$. Furthermore if $\dim M =2$ $h$ is rational.
\end{theorem}
We say that $\beta$ is differentiable in $k$ directions at $h$ if the tangent cone to $\beta$ at $h$ contains a linear space of dimension $k$. We are thus led to ask whether $\beta$ is always differentiable in $k$ directions at a $k$-irrational homology class.  Mather conjectures it is true for $C^{\infty}$ Lagrangians. It is true for all time-periodic Lagrangians on $\T^{1}$ by Theorem \ref{Mather}.  

For autonomous Lagrangians the relevant notion of irrationality is  slightly different. Let $I(h)$ be the dimension of the closure of the image in $\T^{b}$ of $\R h$ instead of $\Z h$. About the relationship between the two notions of irrationality, see Section \ref{irrationality}. Note that  the function $I(h)$  is zero-homogeneous, that is, $I(th)=I(h)$ for all $h \in H_1(M,\R)$ and $t \neq 0$. 

In this context the problem is to show that $\beta$ is always differentiable in $I(h)$ directions at any homology class $h$. This cannot be true in full generality by \cite{BIK}, but it is true for  Finsler metrics on a compact, orientable surface, see \cite{Bangert94} for $M=\T^{2}$ and \cite{ijm}, Corollary 3 for other genera. In \cite{Bangert_Auer} this result is generalized to stable norms in codimension one. 

Note that for autonomous Lagrangians, by \cite{Carneiro}, $\beta$ is always differentiable in the radial direction, that is, for a given $h \in H_1 (M,\R) \setminus \left\{0\right\}$, the function 
\[
\begin{array}{rcl}
\left]0,+\infty \right[ & \longrightarrow & \left]0,+\infty \right[ \\ 
 t & \mapsto & \beta (th)
 \end{array}
 \]
  is $C^{1}$. So the $\beta$ function of an autonomous Lagrangian cannot have vertices except at the zero class. The closest thing to a vertex is a tangent cone which contains no plane, that is, $\beta$ is differentiable in no direction other than radial. 
So far the only positive result for manifolds of dimension greater than two is 
\begin{theorem}[\cite{BIK}]\label{BIK}
The $\beta$-function of a Finsler metric on a manifold $M$ is always differentiable in at least one non-radial direction  at a completely irrational homology class.
\end{theorem}
Note that this theorem contains the result of \cite{Bangert94} because on a two-torus for any non-zero homology class $h$ we have  either $I(h)=1$, in which case $\beta$ is  differentiable in at least one direction, namely the radial one,  the  1-irrational, or $I(h)=2$, that is, $h$ is  completely irrational. 
We extend Theorem \ref{BIK} to 
\begin{theorem}
If the $\beta$ function of an autonomous Lagrangian on a manifold $M$ is not differentiable in any direction other than radial at  $h$, and the energy level that contains the supports of the $(L,h)$-minimizing measures does not meet the zero section of $TM$,  then the support of every $(L,h)$-minimizing measure has Hausdorff dimension $\geq 2I(h)-1$. In particular $2I(h)-1 \leq \dim M$. 
\end{theorem}
Note that this theorem contains Theorem \ref{BIK} because if the Lagrangian is a Finsler metric the only energy level that meets the zero section is the zero level, and it contains only fixed points, whose homology class is zero. Also note that for large enough energy, the corresponding energy level does not meet the zero section.

Here is a layout of the paper. 
 In section \ref{Jacobi} we introduce the main tool of this paper, which is the fibration of $M$ over the circle defined by an integer one-form. In the remaining two sections we prove Theorems 3 and 5 respectively.

\textbf{Acknowledgements} : for many interesting conversations I thank Ivan Babenko,  Gonzalo Contreras and last but not least Osvaldo Osuna who kindly let me use the title of \cite{Osvaldo}. The comments provided by Patrick Bernard and Albert Fathi  led to  vast improvements of  this paper.

\section{Preliminaries}
In this section we briefly recall a few definitions, referring the reader to the bibliography (\cite{Fathi}, and \cite{Bernard_pg} for the time-periodic case) for more information. 
Let $L$ be a complete time-periodic Lagrangian on a closed manifold $M$ of dimension $d$. 
\subsection{The Minimal Action Functional}
Since the $\alpha$-function of $L$ is convex, at every point its graph has a supporting hyperplane. We call face of $\alpha$ the intersection of the graph of $\alpha$ with one of its supporting hyperplane. By Fenchel (a.k.a. convex) duality it is equivalent to study the differentiability of $\beta$  or to study the faces of $\alpha$. 
If $c$ is a cohomology class, we call $F_c$ the largest face of $\alpha$ containing $c$ in its relative interior (it exists by \cite{gafa}), and $V_c$ the underlying vector space of the affine space it generates in $H^{1}(M,\R)$. We call $\tilde{V}_c$ the  vector subspace of $H^{1}(M,\R)\times \R$ generated by  pairs $(c'-c, \alpha (c')-\alpha (c))$ where $c' \in F_c$.
Replacing, if necessary, $L$ by $L-\omega$ where $[\omega]=c$, we only need consider the case when $c=0$.
Likewise, replacing $L$ with $L-\alpha(0)$ we   assume  throughout the paper that $\alpha(0)=0$.
\subsection{Weak KAM preliminaries}
Define, for all $n \in \N$, 
	\[ 
	\begin{array}{rcl} h_n \co \left(M\times \T^1 \right)\times \left(M\times \T^1 \right) & \longrightarrow & \R \\
\left((x,t),(y,s)\right)& \longmapsto & \min \int^{s+n}_{t}L(\gamma,\dot\gamma,t)dt	
\end{array}
\]
where the minimum is taken over all absolutely continuous curves \\
$\gamma \co \left[t,s+n\right]\longrightarrow M$ such that $\gamma (t)=x$ and $\gamma (s+n)= y$. Note that we abuse notation, denoting by the same $t$ an element of $\T^{1}=\R/\Z$ or the corresponding point in $\left[0,1\right[$. The Peierls barrier is then defined as
	\[ \begin{array}{rcl} h \co \left(M\times \T^1 \right)\times \left(M\times \T^1 \right) & \longrightarrow & \R \\
\left((x,t),(y,s)\right)& \longmapsto & \liminf_{n \rightarrow \infty} h_n \left((x,t),(y,s)\right).
\end{array}
\]
The Aubry set is
	\[\Azero := \left\{(x,t) \in M\times \T^1  \co h\left((x,t),(x,t)\right) =0 \right\}.
\]

\begin{definition} 
Let $\tilde{E}_0$ be the set of $(c,\tau) \in H^1 (M\times \T^1,\R)=H^1 (M,\R)\times H^1 ( \T^1,\R)$ such that there exists  a smooth closed one-form $\omega$ on $M\times \T^1$ with $[\omega] = (c,\tau)$ and $\spt (\omega)\cap \Azero = \emptyset$. Let $E_0$ be the canonical projection of $\tilde{E}_0$ to $H^1 (M,\R)$.
\end{definition}
We proved in \cite{soussol} (see \cite{ijm}, Theorem 1 for the autonomous case) the
\begin{theorem}\label{inclusion} The following inclusion holds true : 
	\[E_0 \subset V_0 .
\]
\end{theorem}
\subsection{Irrationality}\label{irrationality}
Recall that $\Gamma$ is the torsion-free part of $H_1 ( M,\Z)$. Let 
\[ J \co H_1 ( M,\R) \longrightarrow H_1 ( M,\R)/\Gamma
\]
denote the canonical projection.
\begin{proposition}\label{ir1}
For any $h \in H_1 ( M,\R)$, the irrationality of $h$ is  the smallest possible dimension of an integer subspace $E_h$ such that $E_h +\Gamma$ contains $h$.
Besides, 
\begin{enumerate}
	\item \label{ir1,1}
	$I(h)$ is  the smallest possible dimension of an integer subspace that contains $h$
	\item \label{ir1,2}
	$I(h)$ is  the dimension of the $\Q$-subspace of $\R$ generated by the coordinates of $h$ in any integer basis of $H_1 ( M,\R)$.
	
\end{enumerate}

\end{proposition}
\proof
Let $h \in H_1 ( M,\R)$ be $k$-irrational with $k\geq 1$.
By the definition of irrationality the closure of $J(\Z h)$ in $H_1 ( M,\R)/\Gamma$ is a finite union of $k$-dimensional tori $G_h$. 
The inverse image of $G_h$ under $J$ is $E_h + \Gamma$, where $E_h$ is a $k$-dimensional vector subspace of $H_1 ( M,\R)$. Now
	\[G_h = \left(E_h + \Gamma\right)/\Gamma \cong E_h /\left( E_h \cap \Gamma \right).
\]
 is compact, so $E_h \cap \Gamma$ is cocompact in $E_h$, that is, $E_h$ is integer. This proves that  $k \geq $ the smallest possible dimension of an integer subspace $E_h$ such that $E_h +\Gamma$ contains $h$.
 
 Conversely, assume $h\in E_h + \Gamma$, with $E_h$ integer of dimension $d$. Then $\Z h \subset  E_h + \Gamma$ whence 
	\[J(\Z h) \subset J(E_h) = E_h /\left( E_h \cap \Gamma \right)
\]
which is a torus of dimension $d$, whence $k\leq d$. This proves that  $k =$ the smallest possible dimension of an integer subspace $E_h$ such that $E_h +\Gamma$ contains $h$.

Let us prove  Statement (\ref{ir1,1}) of the Proposition. Take $h \in H_1 ( M,\R)$. The closure of $J(\R h)$ in $H_1 ( M,\R)/\Gamma$ is an $I(h)$-dimensional tori $G_h$. 
The inverse image of $G_h$ under $J$ is $E_h + \Gamma$, where $E_h$ is an $I(h)$-dimensional vector subspace of $H_1 ( M,\R)$. Now
	\[G_h = \left(E_h + \Gamma\right)/\Gamma \cong E_h /\left( E_h \cap \Gamma \right).
\]
 is compact, so $E_h \cap \Gamma$ is cocompact in $E_h$, that is, $E_h$ is integer. The connected subset $\R h$ contains zero and is contained in $E_h + \Gamma$, so it is contained in the connected component of zero in $E_h + \Gamma$, which is $E_h$.
 This proves that  $I(h) \geq $ the smallest possible dimension of an integer subspace that contains $h$.

Conversely, assume  $h$ is contained in an integer subspace $E_h$ of dimension $k$. Then $E_h$ contains $\R h$ so $J(\R h)$ is contained in the $k$-dimensional torus $E_h /\left( E_h \cap \Gamma \right)$ whence $I(h)\leq k$. This proves that  $I(h)$ is the smallest possible dimension of an integer subspace that contains $h$ .

Now let us prove Statement (\ref{ir1,2}). Take $h \in H_1 ( M,\R)$. Let $a_1,\ldots a_b$ be an integer basis of $H_1 ( M,\R)$, that is, a basis of $H_1 ( M,\R)$ all of whose elements are in $\Gamma$ and set $h=\sum^{b}_{i=1}\lambda_i a_i$. 
Let 
	\[k:= \dim \mbox{Vect}_{\Q}\left\{\lambda_1,\ldots \lambda_b\right\}.
\]
Modulo a permutation of $\lambda_1,\ldots \lambda_b$ we may assume that 
\begin{eqnarray*}
	\lambda_{k+1}&=&\sum^{k}_{j=1}r_{j,k+1}\lambda_j\\
	\vdots & \vdots & \vdots \\
	\lambda_{b}&=&\sum^{k}_{j=1}r_{j,b}\lambda_j
\end{eqnarray*}
where all the $r_{j,i}$ are rationals. Let $N$ be their least  common multiple. Then 
\begin{eqnarray*}
	Nh &=& \sum^{b}_{i=1}N\lambda_i a_i \\
	&=& \sum^{k}_{i=1}N\lambda_i a_i +\sum^{b}_{i=k+1}\left(\sum^{k}_{j=1}N r_{j,i}\lambda_j\right) a_i \\
	&=& \sum^{k}_{i=1}\lambda_i\left[N a_i  +   \sum^{b}_{l=k+1}N r_{i,l} a_l  \right]
\end{eqnarray*}
that is, $h$ lies in the integer subspace of $H_1 ( M,\R)$ generated by the $k$ classes
	\[N a_i  +   \sum^{b}_{l=k+1}N r_{i,l} a_l ,\  i=1,\ldots k
\]
so $I(h) \leq k$. 

Conversely, let $E_h$ be an integer suspace of $H_1 ( M,\R)$ of dimension $I(h)$ containing $h$.
Let $b_1,\ldots b_d$ be an integer basis of $E_h$ and set $h=\sum^{I(h)}_{j=1}\nu_j b_j$. 
Since $b_1,\ldots b_d$ are elements of $\Gamma$ there exist rational numbers $r_{i,j}$, $i=1,\ldots b$,  $j=1,\ldots I(h)$ such that 
	\[b_j = \sum^{b}_{i=1}r_{i,j}a_i,\ \forall j=1,\ldots I(h).
\]
Then 
	\[h =\sum^{b}_{i=1}\left(\sum^{j=1}_{I(h))}r_{i,j}\nu_j \right)a_i
\]
so the coordinates of $h$ in the basis $a_1,\ldots a_b$ are linear combinations, with rational coefficients, of the $\nu_j$, $j=1,\ldots I(h)$, thus the $\Q$-subspace of $\R$ they generate has dimension at most $I(h)$.
\qed
\section{Integrality of $\tilde{V}_0$}
We begin by explaining why the integrality of  $\tilde{V}_0$ matters in the problem we study. 
\subsection{Consequence of the integrality of  $\tilde{V}_0$}\label{tilde_entier}
Let us assume that $\tilde{V}_0$ is an integer subspace of $H^1 (M\times \T^1,\R)$. Take cohomology classes
$c_1,\ldots c_k$ in $F_0$ and $\omega_1,\ldots \omega_k$ closed one-forms on $M$ such that 
\begin{itemize}
	\item 
	$[\omega_i]=c_i, i=1\ldots k$
	\item
	$\forall i=1\ldots k, \exists n_i \in \Z, n_i (c_i, \alpha (c_i)) \in H^1 (M\times \T^1,\Z)$,
	in particular $n_i  \alpha (c_i) \in \Z$
	\item
	the classes $(c_i, \alpha (c_i)), i=1\ldots k$ form a basis of $\tilde{V}_0$.
\end{itemize}
Complete $n_i (c_i, \alpha (c_i)),i=1\ldots k$ to an integer basis $B$ of $H^1 (M\times \T^1,\Z)$.
Let $\mu$ be a minimizing measure, for $L$, then since $c_i \in F_0, i=1\ldots k$, $\mu$ is also $L+\omega_i$-minimizing. Thus for $i=1,\ldots k$ we have
\begin{eqnarray*}
\int_{TM\times \T^1} \left(L+\omega_i \right)d\mu &=& \alpha (c_i)\\
\int_{TM\times \T^1} Ld\mu &=& 0
\end{eqnarray*}
(recall that we work under the assumption that $\alpha (0)=0$) whence
	\begin{equation}\label{integrality1}\int_{TM\times \T^1} \omega_i d\mu = \alpha (c_i)	\in \Q
\end{equation}
but on the other hand
	\begin{equation}\label{integrality2}\int_{TM\times \T^1} \omega_i d\mu =<c_i,\left[\mu \right]>
\end{equation}
so the first $k$ coordinates of $\left[\mu \right]$ in the basis $h_1,\ldots h_b$ of $H_1 (M,\R)$ dual to $B$ are rational. So for some $N \in \Z$, $Nh$ belongs to $\Gamma +\mbox{Vect}(h_{k+1},\ldots h_b )$. Since $\mbox{Vect}(h_{k+1},\ldots h_b )$ is integer and of dimension $b-k$, the irrationality  of $ \left[\mu \right]$ is at most $b-k$. 

By transposition we see that if $h \in H_1 (M,\R)$ is $k$-irrational, and if $c \in H^{1} (M,\R)$ is a subderivative to $\beta$ at $h$, with $\tilde{V}_c$ integer, then 
	\[\dim F_c  \leq \dim \tilde{V_c} \leq b-k
\]
whence $\beta$ is differentiable at $h$ in at least $k$ directions.

\subsection{Further consequence of the integrality}\label{Jacobi}
Assume $\tilde{V}_{0}$ contains an integer point. Let $c \in F_0$ be such that, for some $\lambda \in \R$ 
	\[ \lambda \left(c, \alpha(c) \right) \in H_1 \left(M \times \T^{1}, \Z \right).
\]
Let $\omega$ be a smooth one-form on $M$ such that $[\omega]=c$. Fix once and for all an origin $(x_0,0) \in \Azero \cap \left\{t=0\right\}$ and consider the map :
\[
\begin{array}{rcl}
\phi_{\omega} \co M \times \T & \longrightarrow & \R/\frac{1}{\lambda}\Z  \\
(x,t) & \longmapsto & \int^{(x,t)}_{(x_0,0)} \left[\omega -\alpha(c)dt \right]	
\end{array}
\]
 where the integral is taken along any path in $M\times \T $ joining $(x_0,0)$ to $(x,t)$.
 
 To emphasize the dependance on the Lagrangian we write $h(L)$ or $h(L-\omega)$ the corresponding Peierls barriers.

For brevity we denote
\begin{eqnarray*}
u_1 (x,t) & :=  &	h(L-\omega)\left( (x_0,0),(x,t) \right)\\
u_0 (x,t) & :=  &	h(L)\left( (x_0,0),(x,t) \right).
\end{eqnarray*}
Define 
\[
\begin{array}{rcl}
\Phi_{\omega} \co M \times \T & \longrightarrow & \R/\frac{1}{\lambda}\Z  \\
(x,t) & \longmapsto & 	\phi_{\omega}(x,t)-u_0 (x,t)+u_1 (x,t).
\end{array}
\]
We have 
\begin{proposition}\label{constante_classes_statiques}
The map $\Phi_{\omega}$ is constant on the static classes of $L$.
\end{proposition}
\proof
First we need a
\begin{lemma} \label{memes_classes_statiques}
The static classes for $L-\omega$ are the same as the static classes for $L$, and furthermore the curves realizing the liminf are the same.
\end{lemma}
\begin{remark}
By \cite{ijm}, Proposition 6 (see \cite{Bernard_Fourier} for the time-periodic case) the Aubry sets for $L$ and $L-\omega$ are the same.
\end{remark}
Take $(x_1,t_1)$ and $(x_2,t_2)$ in the same static class for $L$. We have 
	\[h(L)\left( (x_1,t_1),(x_2,t_2) \right)+h(L)\left( (x_2,t_2),(x_1,t_1)\right)=0
\]
that is to say, for all $n \in \N$ there exist  two increasing sequence of integers $k_n, l_n,\  n \in \N$ and curves 
\[
\begin{array}{rcl}
\gamma_n  \co \left[t_1,t_2+k_n\right]& \longrightarrow & M \\
\delta_n  \co \left[t_2,t_1+l_n\right]& \longrightarrow & M
\end{array}
\]
 such that
$\gamma_n (t_1)=\delta_n (t_1+l_n) =x_1$, $\gamma_n (t_2+k_n)=\delta_n (t_2) =x_2$ and 
	\[\int^{t_2+k_n}_{t_1}L(\gamma_n, \dot{\gamma}_n,t)dt +\int^{t_1+l_n}_{t_2}L(\delta_n, \dot{\delta}_n,t)dt\longrightarrow 0
\]
when $n \rightarrow \infty$.
Since $c \in F_0$, and $F_0$ contains $0$ in its interior, we can find $\theta_1 >0$ such that 
$c_1 := -\theta_1 c \in F_0$ so, denoting $\theta := (1+\theta_1)\inv$, we have $(1-\theta)c + \theta c_1 =0$.
Thus, $\alpha$ being affine in $F_0$, we have 
	\[ 0= \alpha (0)= \alpha ((1-\theta)c)+\alpha (\theta c_1 )=(1-\theta)\alpha (c)+\theta \alpha (c_1 ).
\]

Now if $v$, $w$ are weak KAM solutions for $L-\omega$ and $L- \omega_1$ respectively (with 
$\omega_1 := -\theta_1 \omega$), we have
\[
 \int^{t_2+k_n}_{t_1}\left(L-\omega +\alpha ((c) \right)(\gamma_n, \dot{\gamma}_n,t)dt +\int^{t_1+l_n}_{t_2}\left(L-\omega +\alpha (c) \right)(\delta_n, \dot{\delta}_n,t)dt 
 \]
\[
 \geq  v(x_2,t_2)-v(x_1,t_1)+	v(x_1,t_1)-v(x_2,t_2)=0\]
\[ 
\int^{t_2+k_n}_{t_1}\left(L-\omega_1 +\alpha ( c_1) \right)(\gamma_n, \dot{\gamma}_n,t)dt +\int^{t_1+l_n}_{t_2}\left(L-\omega_1 +\alpha ( c_1) \right)(\delta_n, \dot{\delta}_n,t)dt 
\]
\[
\geq  w(x_2,t_2)-w(x_1,t_1)+	w(x_1,t_1)-w(x_2,t_2)=0.
\]
Summing $\theta$ times the first inequality with $1-\theta$ times the second we get a sum of two non-negative terms which converges to zero, so each term must converge to zero. Consequently
\[h(L-\omega)\left( (x_1,t_1),(x_2,t_2) \right)+h(L-\omega)\left( (x_2,t_2),(x_1,t_1)\right)=0
\]
which proves that $(x_1,t_1)$ and $(x_2,t_2)$ are in the same static class for $L-\omega$. Lemma \ref{memes_classes_statiques} follows by swapping $L$ and $L-\omega$.

\qed

Now let us prove Proposition \ref{constante_classes_statiques}. We have, re-using the notations of Lemma \ref{memes_classes_statiques},
\begin{eqnarray*}
h(L)\left( (x_1,t_1),(x_2,t_2) \right) &=& \lim_{n \rightarrow \infty}\int^{t_2+k_n}_{t_1}L (\gamma_n, \dot{\gamma}_n,t)dt\\ 	
h(L-\omega)\left( (x_1,t_1),(x_2,t_2) \right)&=& \lim_{n \rightarrow \infty}\int^{t_2+k_n}_{t_1}\left(L-\omega +\alpha (c) \right)(\gamma_n, \dot{\gamma}_n,t)dt 	
\end{eqnarray*}
whence by substraction
\begin{eqnarray*}
\lim_{n \rightarrow \infty}\int^{t_2+k_n}_{t_1}\left[\omega (\gamma_n, \dot{\gamma}_n,t)-\alpha(c)\right]dt &=&
h(L)\left( (x_1,t_1),(x_2,t_2) \right)\\
&&-	h(L-\omega)\left( (x_1,t_1),(x_2,t_2) \right).	
\end{eqnarray*}

Thus
	\[\phi_{\omega} (x_2,t_2)-\phi_{\omega} (x_1,t_1)=
	\]
	\[
	h(L)\left( (x_1,t_1),(x_2,t_2) \right)-	h(L-\omega)\left( (x_1,t_1),(x_2,t_2) \right)\mbox{ mod } \frac{1}{\lambda}\Z.
\]
On the other hand, since $(x_1,t_1)$ and $(x_2,t_2)$ are in the same static class, we have (cf \cite{Bernard_pg}, Lemma 4.1)
\begin{eqnarray*}
h(L)\left( (x_0,0),(x_2,t_2) \right)&=&h(L)\left( (x_0,0),(x_1,t_1) \right)+h(L)\left( (x_1,t_1),(x_2,t_2) \right)\\
h(L-\omega)\left( (x_0,0),(x_2,t_2) \right)&=&h(L-\omega)\left( (x_0,0),(x_1,t_1) \right)\\
&&+h(L-\omega)\left( (x_1,t_1),(x_2,t_2) \right)
\end{eqnarray*}
that is,
\begin{eqnarray*}
h(L)\left( (x_1,t_1),(x_2,t_2) \right) &=& u_0(x_2,t_2)-u_0 (x_1,t_1) \\
h(L-\omega)\left( (x_1,t_1),(x_2,t_2) \right) &=& u_1(x_2,t_2)-u_1 (x_1,t_1)	
\end{eqnarray*}
and finally 
	\[\phi_{\omega} (x_2,t_2)-u_0(x_2,t_2)+u_1 (x_2,t_2)=
	\]
	\[\phi_{\omega} (x_1,t_1)-u_0(x_1,t_1)+u_1 (x_1,t_1) \mbox{ mod }\frac{1}{\lambda}\Z.
\]
\qed

\begin{proposition}\label{Holder}
The map $\Phi_{\omega}$ satisfies a   H\H{o}lder condition of order two along $\Azero $.
 \end{proposition}
\proof
Since $c\in F_0$, we have $\Azero (L-\omega)= \Azero (L)$. Besides, the maps
\[
\begin{array}{rcl}
 M \times \T^{1} & \longrightarrow & \R \\
(x,t) & \longmapsto & h(L-\omega)\left( (x_0,0), (x,t) \right)\\
(x,t) & \longmapsto &	h(L)\left( (x_0,0),(x,t) \right)	
\end{array}
\]
are weak KAM solutions, so at every point $(x,t)$ of $\Azero$ there exists a $v \in T_xM$ such that 
\begin{eqnarray*}
\frac{\partial}{\partial x}h(L-\omega)\left( (x_0,0),(x,t) \right)&=&
\frac{\partial L}{\partial v}	(x,v,t) - \omega_x \\
\frac{\partial}{\partial x}h(L)\left( (x_0,0),(x,t) \right)&=&
\frac{\partial L}{\partial v}	(x,v,t).
\end{eqnarray*}
Furthermore (see \cite{soussol}) 
	\[\frac{\partial}{\partial t}h(L)\left( (x_0,0),(x,t)\right) - \frac{\partial}{\partial t}h(L-\omega)\left( (x_0,0),(x,t) \right)= -\alpha(c)dt.
\]
Let $(U,f)$ be a coordinate chart for $M\times \T$. 
We have (see \cite{Fathi}, 4.5.5 for the autonomous case)
\begin{lemma}
There exists $K$ such that for any weak KAM solution $u$ for $L$, for any $X$ in $U$ such that $f(X) \in \Azero$,  for any $Y$ in $U$, we have
	\[\left|u(f(Y))-u(f(X)) -d_{f(X)}u.d_X f(Y-X)\right|\leq K\left\| Y-X  \right\|^{2}
\]
\end{lemma}
\proof
Take $X,Y$ in $U$ such that $f(X) \in \Azero $. 
Let $u_-$ be a backward weak KAM solution for $L$, and let $u_+$ be its conjuguate, so  $u_+$ is a forward weak KAM solution for $L$. By \cite{Bernard_pg} there exists a real number $K$ such that all backward (resp. forward) weak KAM solutions are K-semi-concave (resp. semi-convex).  Since $f(X) \in \Azero $, $u_-$ and 
$u_+$ are differentiable, and have the same derivative,  at $f(X)$. So by semi-concavity (resp. semi-convexity)
\begin{eqnarray*}
u_-(f(Y))-u_-(f(X))&\leq & d_{f(X)}u.d_{X}f(Y-X)	+K\left\| Y-X  \right\|^{2}\\
u_+(f(Y))-u_+(f(X))&\geq & d_{f(X)}u.d_{X}f(Y-X)	+K\left\| Y-X  \right\|^{2}.
\end{eqnarray*}
On the other hand, $u_-$ and $u_+$ being weak KAM conjuguate, $u_-(f(X))=u_+(f(X))$ because $f(X) \in \Azero $ and $u_-(f(Y)) \geq u_+(f((Y))$ so 
	\[u_-(f(Y))-u_-(f(X)) \geq u_+(f(Y))-u_+(f(X))
\]
which combines with the former two inequalities to prove the lemma.

\qed

Now let us prove Proposition \ref{Holder}. By our previous lemma we have, for any $X,Y$ such that $f(X),f(Y) \in \Azero$,
	\[	\left| u_1 (f(Y))-u_0 (f(Y)) - u_1 (f(X))+u_0 (f(X))	+\left(\omega_{f(X)}-\alpha(c)dt \right).d_{X}f (Y-X)     \right| 
\]
	\[\leq  2K\left\|Y-X\right\|^{2}.
\]
On the other hand, since $\omega$ is smooth, by Taylor's formula we have for some $K'$ :
	\[ \left| \phi_{\omega}(f(Y))-\phi_{\omega}(f(X))-
	\left[\omega_{f(X)}-\alpha(c)dt \right].d_{X}f (Y-X)    \right|\leq  K'\left\|Y-X\right\|^{2}
\]
thus
	\[ \left| \Phi_{\omega}(f(Y))-\Phi_{\omega}(f(X))    \right|\leq  (2K+K')\left\|Y-X\right\|^{2}
\]
which proves the Proposition.
\qed

\subsection{}
In this paragraph we prove a partial converse to Theorem \ref{inclusion}.
Consider the equivalence relation on $M\times \T$ defined by $(x,t)\approx (y,s)$ if and only if $h((x,t),(y,s))+h((y,s),(x,t))=0$. The quotient $\left(M\times \T\right)/\approx$ is a metric space with distance $d(\overline{(x,t)},\overline{(y,s)})=h((x,t),(y,s))+h((y,s),(x,t))$, where $\overline{(x,t)}$ is the equivalence class of $(x,t)$. The image of $\Azero$ in $\left(M\times \T\right)/\approx$  is the quotient Aubry set $A_0$ of \cite{Mather02}. 

\begin{proposition}\label{integer_countable}
If the quotient Aubry set $A_0$ is totally disconnected, and $\tilde{V_0}$ contains an integer point $(c,\tau)$,  then  
\[(c,\tau) \in \tilde{E}_0 .
\] 
In particular, if $A_0$ is totally disconnected, and $\tilde{V_0}$ is  integer, then 
\[ E_0 = V_0 .
\] 
\end{proposition}
\begin{remark}
By \cite{Mane96} for a generic Lagrangian there exists a unique minimizing measure, and in particular a unique static class, so $A_0$ is a point. Also, by 	\cite{Mather02}, if $\dim M=2$, $A_0$ is always totally disconnected. 
 \end{remark}
 \proof
 Take $c$ in $F_0$, such that some multiple of $\left(c,\alpha(c)\right)$  is integer. Take a one-form $\omega$ such that 
 $\left[\omega\right]=c$ and define $\Phi_{\omega}$ as above. 
 
 First assume that 
	\[\Phi_{\omega} (\Azero ) \neq \T^{1}.
\]
Since $\Azero$ is compact, so is 
 $\Phi_{\omega}\left(\Azero \right)$. Therefore 
$\Phi_{\omega}\left(\Azero \right)$ misses an open set $V_1$. Then we can find a map $\Phi$ which is $C^{1}$, homotopic to   $\Phi_{\omega}$ and misses an open set $V\subset V_1$.
Take a closed one form $\eta$  
generating the integral cohomology of the circle $\R/\lambda\inv\Z$, supported
in $V$. The pull-back $\Phi^* (\eta)$ of this form
by $\Phi$  defines a closed one-form in $M\times \T$, supported away from
$\Azero $. Since $\eta$ is cohomologous to the constant one-form $dx$ on the circle, we have
	\[ \left[\Phi^* (\eta) \right]=\left[\Phi^* (dx) \right]
\]
but $\Phi^* (dx)$ is cohomologous to $\Phi^{\ast}_{\omega}(dx)$ since $\Phi$  is homotopic to   $\Phi_{\omega}$. By the definition of $\Phi_{\omega}$, $\Phi_{\omega}^* (dx)$ is 
cohomologous to $(\omega,\alpha(c)dt)$,  which implies $c \in E_0$. 

So what remains to do now is assume that $\Phi_{\omega}$ restricted to $\Azero$ is surjective and use the hypothesis of total disconnectedness to modify  $\Phi_{\omega}$. 
\begin{remark}
By \cite{Rifford}, Lemma 2, if $\dim M=2$, $\Phi_{\omega}\left(\Azero \right)$ has Lebesgue measure zero, so $\Phi_{\omega}$ restricted to $\Azero$ is not surjective.
\end{remark}
Assume, up to a translation, 
that   $ \Phi_{\omega}(x_0,0)=0$. We are looking for a map $\Phi$ which  has the same effect as $\Phi_{\omega}$ on cohomology classes and misses a  neighborhood  of the point $(2\lambda)\inv$ in the circle.

By \cite{HY}, Theorem 2.15,  in a totally disconnected metric space, every point has a basis of open-closed neighborhoods. Thus, since  the canonical  projection  from $M$ to $\left(M\times \T\right)/\approx$ is continuous, for every static class $S$ contained in $\Phi _{\omega}\inv (1/2 \lambda)\cap \Azero$ there exists an open (in $\Azero$) and closed neighborhood $V_S$ of $S$ such that $\Phi_{\omega}(V_S) \subset \left[1/4\lambda,3/4\lambda\right]$. Cover the compact set $\Phi _{\omega}\inv (1/2 \lambda)\cap A_0$ with the $V_S$'s and take a finite subcover $V_{1},\ldots V_{n}$. Then the reunion $V_{1}\cup \ldots \cup V_{n}$
is an open (in $\Azero$) and closed neighborhood of $\Phi _{\omega}\inv (1/2 \lambda)\cap \Azero$. So its complement in $\Azero$ is closed. This shows that $\Azero$ splits as disjoint union of two compact subsets $K_1$ and $K_2$ such that  
\begin{itemize}
	\item $K_1 \supset\Phi _{\omega}\inv (1/2 \lambda)\cap \Azero$,
	\item $\Phi _{\omega}(K_1) \subset \left[1/4\lambda,3/4\lambda\right]$ and
	\item $1/2\lambda \notin \Phi _{\omega}(K_2) $. 
\end{itemize}
  
Thus there exists a $C^{1}$ function $f$ on $M\times \T$
such that $f$ is zero on $K_2$ and one on $K_1$.

Replace $\omega$ with $\omega -(2\lambda)\inv df$ in the definition of $\Phi_{\omega}$ (without changing $u_1$), and call $\Phi$ the resulting map. Then, if $x \in K_2$, $\Phi (x)=\Phi_{\omega}(x) \neq 1/2\lambda$, and if $x \in K_1$, 
	\[\Phi (x)=\Phi_{\omega}(x) -\frac{1}{2\lambda}f(x)= \Phi_{\omega}(x)-\frac{1}{2\lambda} \neq \frac{1}{2\lambda}
\]
 which shows that $\Phi$ restricted to $K$ is not surjective and brings us back to the first case.
 This proves Proposition \ref{integer_countable}.
\qed

 \subsubsection{}
 Let us now consider the case when there are so many static classes that 	
 \[\Phi_{\omega} \left( \Azero  \right)= \R/\Z.
\]
 Note that, since $\Phi_{\omega}$ is constant on the static classes,  a fortiori it is constant on the orbits of points in $\Azero$, so, taking the initial point of each orbit in $\Azero$,  
	\[\Phi_{\omega} \left( \Azero \cap \{t=0\} \right)= \R/\Z.
\]
 
  The classical lemma below then shows that $\Azero \cap \{t=0\} \subset M\times \left\{0\right\}$ has Hausdorff dimension $\geq 2$. 

\begin{lemma} \label{dimHausdorff}
Let $A$ be a subset of $\R^n$ and let $C, \alpha$ be positive real numbers.
Take $f \co A \longrightarrow \R$  a function such that 
$\forall x,y \in A$, we have $|f(x)-f(y)| \leq C |x-y|^{\alpha}$.
Then, denoting by $\mathcal{H}^{s}$ the Hausdorff measure of dimension $s$, we have 
\[
\mathcal{H}^{\frac{s}{\alpha}}(f(A)) \leq 
C^{\frac{s}{\alpha}} \mathcal{H}^{s}(A).
\]
\end{lemma}
\proof
If $(A_i)_i$ is a $\delta$-covering of $A$, then $(f(A_i))_i$ is a 
$C\delta ^{\alpha}$-covering of $f(A)$, and 
\[
\sum_{i=1}^{\infty} \mbox{diameter}(f(A_i))^{\frac{s}{\alpha}} \leq
\sum_{i=1}^{\infty} C^{\frac{s}{\alpha}} \mbox{diameter}(A_i)^s .
\]
Taking the infimum over all $\delta$-coverings of $A$, and letting $\delta$
go to zero, we get the required conclusion.
\qed
 \section{Proof of Theorem 3}
 Let us now assume that $V_0$ is integer whithout specification on $\tilde{V}_0$.
 Take cohomology classes
$c_1,\ldots c_k$ in $F_0$ and $\omega_1,\ldots \omega_k$ smooth closed one-forms on $M$ such that 
\begin{itemize}
	\item 
	$[\omega_i]=c_i, i=1\ldots k$
	\item
	$\forall i=1\ldots k, \exists n_i \in \Z, n_i c_i \in H^1 (M,\Z)$,
	\item
	the classes $c_i, i=1\ldots k$ form a basis of $V_0$.
\end{itemize}
Define, for each $i$,
\[
\begin{array}{rcl}
\phi_{i} \co M  & \longrightarrow & \R/\frac{1}{n_i}\Z  \\
x & \longmapsto & \int^{x}_{x_0} \omega_i 	
\end{array}
\]
 where the integral is taken along any path in $M$ joining $x_0$ to $x$.
  For brevity we denote
\begin{eqnarray*}
u_{1,i} (x) & :=  &	h(L-\omega_i)\left( (x_0,0),(x,0) \right)\\
u_0 (x) & :=  &	h(L)\left( (x_0,0),(x,0) \right).
\end{eqnarray*}
Define 
\[
\begin{array}{rcl}
\Phi_{i} \co M  & \longrightarrow & \R/\frac{1}{n_i}\Z  \\
x & \longmapsto & 	\phi_{i}(x)-u_0 (x)+u_{1,i} (x).
\end{array}
\]
 Consider the map 
\[
\begin{array}{rcl}
\Phi \co M & \longrightarrow & \T^{k} =\bigoplus^{k}_{i=1}\R/\frac{1}{n_i}\Z \\
x & \longmapsto & \left( \Phi_{i} 	\right)_{i=1,\ldots k} 
\end{array}
\]
 It can be proved as in Proposition \ref{Holder} that $\Phi$ satisfies a   H\H{o}lder condition of order two along $\Azero \cap \{t=0\}$.
 The difference with the previous section is that the vector $(\alpha(c_i))_{i=1,\ldots k}$ may be irrational so $\Phi $ may not be constant on the static classes.
 Take $\gamma \co \R \longrightarrow M$ an extremal  such that $(\gamma (t), \dot{\gamma} (t),t)$ is contained in $\Azero$ for all $t$. We have, for all $n$ in $\Z$,
 \begin{eqnarray*}
 \int^{n}_{t=0}\left(L-\omega_i +\alpha(c_i)\right)(\gamma (t), \dot{\gamma} (t),t)dt 
 &=& u_{1,i}(\gamma(n))-u_{1,i}(\gamma(0))\\
 \int^{n}_{t=0}\left(L\right)\left(\gamma (t), \dot{\gamma} (t),t\right)dt &=& u_{0}(\gamma(n))-u_{0}(\gamma(0))
\end{eqnarray*} 
 Substracting the last two equations we get 
	\[\Phi \left(\gamma (n) \right)-\Phi \left(\gamma (0) \right) = 
	n\left(\alpha(c_i) \mbox{ mod }\frac{1}{n_i}\Z \right)_{i=1,\ldots k}.
\]

 Now if the vector $(\alpha(c_i))_{i=1,\ldots k}$ is $p$-irrational the set 
 
	\[\left\{n(\alpha(c_i))_{i=1,\ldots k} \co n \in \Z \right\}
\]
 is dense in a subtorus of dimension $p$. Thus by the  H\H{o}lder property the closure of $ \gamma \left(\Z \right)$ in $M$ has dimension $\geq 2p$, whence the dimension of $M$ itself is at least $2p$ and the support of any minimizing measure, which is a subset of $TM\times \T^{1}$, has dimension $\geq 2p+1$. 
\subsection{Vertices}
Suppose $\beta$ has a vertex at some homology class $h$. Then there exists a cohomology class $c$ such that $c$ is a subderivative to $\beta$ at $h$ and $ V_c = H^{1}(M,\R)$. In particular $V_c$ is an integer subspace so this situation is contained in the case we just considered.  Pick an $h$-minimizing measure $\mu$ and an integer basis $c_1,\ldots c_b$ of $H^{1}(M,\R)$. Let $h_1,\ldots h_b$ be the basis  of $H_{1}(M,\R)$ dual to $c_1,\ldots c_b$.  The coordinates of $h$ in the basis $h_1,\ldots h_b$ are $<c_i,h>,\;i=1,\ldots b$. Note that since $c_1,\ldots c_b$ is an integer basis, some mutiple of its dual basis is an integer basis of $H_{1}(M,\R)$, so the irrationality of $h$ is just that of the vector $<c_i,h>,\;i=1,\ldots b$. Now,  $<c_i,h>-\alpha(c_i)=0$ for all $i=1,\ldots b$ by Equations (\ref{integrality1}) and (\ref{integrality2}). So the irrationality $p$ of $h$ is exactly that of $(\alpha(c_i))_{i=1,\ldots k}$, which is $\leq 1/2 \dim M$. Besides, the support of $\mu$ has Hausdorff dimension $\geq 2p+1$. This finishes the proof of Theorem \ref{vertices}, but for the two-dimensional case.

\subsection{The two-dimensional case}
We shall prove that if the dimension of $M$ is two and $V_0$ is integer, then so is $\tilde{V}_0$. More specifically we prove that if $c \in V_0$ is integer, then so is some multiple of $(c,\alpha(c))\in \tilde{V}_0$. As a bonus, using Proposition \ref{integer_countable}, we get the following
\begin{corollary}
If the dimension of $M$ is two and $V_0$ is integer, then $E_0=V_0$.
\end{corollary}
Let $c \in F_0$ be such that for some $n \in \Z$, $nc \in H^{1}(M,\Z)$. Let $\omega$ be a smooth closed one-form on $M$ with cohomology $c$. Denote
\begin{eqnarray*}
u_1 (x) & :=  &	h(L-\omega)\left( (x_0,0),(x,0) \right)\\
u_0 (x) & :=  &	h(L)\left( (x_0,0),(x,0) \right).
\end{eqnarray*}
Define 
\[
\begin{array}{rcl}
\phi_{\omega} \co M  & \longrightarrow & \R/\frac{1}{n}\Z  \\
x & \longmapsto & \int^{x}_{x_0} \omega 	
\end{array}
\]
 where the integral is taken along any path in $M$ joining $x_0$ to $x$, and 
\[
\begin{array}{rcl}
\Phi_{\omega} \co M  & \longrightarrow & \R/\frac{1}{n}\Z  \\
(x,t) & \longmapsto & 	\phi_{\omega}(x,t)-u_0 (x)+u_1 (x).
\end{array}
\]
Then $\Phi_{\omega}$ is $2$-H\H{o}lder along $\Azero$, and if $\alpha(c)$ is irrational, then $\Phi_{\omega}(\Azero)$ is surjective, which contradicts Lemma 2 of \cite{Rifford}. Therefore $\alpha(c)$ is rational. In particular if $\beta$ has a vertex at some $h$, then 
$<c_i,h>$ is rational for any integer basis $c_i,\; i=1,\ldots , b$ of $H^{1}(M,\R)$. Thus $h$ is rational. This proves the two-dimensonal case of Theorem \ref{vertices}.
 \qed
 \section{Proof of Theorem 5}
 Let $h \in H_1 (M,\R)$ be such that the tangent cone to $\beta$ at $h$ contains no plane. Let $c \in H^1 (M,\R)$ be such that $<c,h>=\alpha (c) + \beta (h)$, that is, $c$ defines a supporting hyperplane to $\beta$ at $h$. Modifying $L$ by a closed one-form if $c \neq 0$, we may assume $c=0$. Then $V_0$ has codimension one or zero. 
 \begin{lemma}
 If the codimension of $V_0$ is zero, then $h=0$.
 \end{lemma}
 \proof  Since $0 \in H^{1}(M,\R)$ is a subderivative to $\beta$ at $h$, and $L$ is autonomous,  by \cite{Carneiro} the support of any $h$-minimizing measure is contained in the energy level $0=\alpha(0)$. Since 
 any $c \in F_0$ is also a subderivative to $\beta$ at $h$, the support of any $h$-minimizing measure is contained in the energy level $\alpha(c)$ so we have 
	\[0=\alpha(c )\; \forall c \in F_0.
\]
In other words, the faces of $\alpha$ are contained in the level sets of $\alpha$ for an autonomous Lagrangian.
Besides, since $c$ and $0$ are subderivatives to $\beta$ at $h$, we have 
\begin{eqnarray*}
<c,h>&=& \alpha (c) + \beta(h)=\beta(h)	\\
<0,h>&=&   \beta(h)	
\end{eqnarray*}
(recall that $\alpha (0)=0$) whence 
	\[<c,h> =<0,h>=0 \; \forall c \in F_0.
\]
  Therefore, if $ V_0 = H^{1}(M,\R)$ we have $<c,h> =0 \; \forall c \in H^{1}(M,\R)$ so $h=0$. 
 \qed

 Since $I(0)=0$, there is nothing to prove in that case so  we may assume $h\neq 0$ and $V_0$ has codimension one. 
 
 We would like an ergodic $h$-minimizing measure. Such a measure need not exist because $(h,\beta(h))$ may not be an extremal point of the epigraph of $\beta$, that is, $h$ may  lie in the relative interior of a face of $\beta$. Such a face must be radial, i.e. contained in $\left\{ th \co t\in \right[ 0,+\infty \left[ \right\}$, for $\beta$ is not differentiable at $h$ in any direction but the radial one. Then the face has at least one non-zero extremal point $th$ with $t\in \left] 0,+\infty \right[$. Since $I(th)=I(h)$, we may, for the purpose of proving Theorem 5, assume that $h$ itself is an extremal point of $\beta$. Then by \cite{Mane92} there exists an ergodic $h$-minimizing measure $\mu$.
 
 Let $c_1,\ldots c_{b-1}$ be a basis of $V_0$ and let $\omega_1,\ldots \omega_{b-1}$ be smooth closed one-forms on $M$ sucht that $[\omega_i]=c_i, i=1\ldots k$. 
 Consider the functions
\begin{eqnarray*}
u_{1,i} (x) & :=  &	h(L-\omega_i)\left( x_0,x \right)\\
u_{0} (x) & :=  &	h(L)\left( x_0,x \right).
\end{eqnarray*}	
for $1,\ldots b-1$, and replace each one-form  $\omega_i$ by 
	\[ \omega'_i := \omega_i + du_{1,i} - du_{0}.
\]
This is an almost everywhere defined, integrable one-form.
 Pick $a_1,\ldots a_{b}$ an integer basis of $H^{1}(M,\R)$ and a matrix 
	\[\Lambda = \left(\lambda_{ij}\right)\in M_{b,b-1}\left(\R \right) \mbox{ such that }
	c_j=\sum^{b}_{1}\lambda_{ij}a_i \; \forall j=1, \ldots b-1
	\]
then the rank of $\Lambda$ is $b-1$ so there exist representatives  $\eta_1,\ldots \eta_{b}$ of $a_1,\ldots a_{b}$ such that 
	\[ \left(\omega'_1,\ldots \omega'_{b-1}\right)=\Lambda \left(\eta_1,\ldots \eta_{b}\right)
\]
 
 Consider the map 
 \[
\begin{array}{rcl}
\Phi \co M & \longrightarrow & \T^{b}  \\
x & \longmapsto & \left( \int^{x}_{x_0} \eta_i 	\right)_{i=1,\ldots b} \mbox{ mod}\Z.
\end{array}
\]
 Pick an orbit $\gamma \co  \R \longrightarrow M$ contained in the support $\spt \mu$ of $\mu$. We have 
	\[ d\Phi (\dot{\gamma}(t))=\left(\eta_i (\dot{\gamma}(t))  \right)_{i=1,\ldots b}
	\mbox{ so } \Lambda \left(d\Phi (\dot{\gamma}(t))  \right)=0
\]
 and $\Phi (\gamma (\R))$ is a connected subset of  a straight line in $\T^{b}$. We claim that for $\mu$-almost every orbit $\gamma$ the closure of 
 $\Phi (\gamma (\R))$ in $\T^{b}$ is a subtorus of dimension  $I(h)$.
 Let $\overline{M}$ be the Abelian cover  of $M$, that is, the cover whose group of Deck transformations is isomorphic to $\Gamma$.
 We denote by $h$ the Deck transformation corresponding to the homology class $h$.
 Lift the map $\Phi$ to  $\overline{M}$  :
\[
\begin{array}{rcl}
\overline{\Phi} \co \overline{M} & \longrightarrow & \R^{b}  \\
x & \longmapsto & \left( \int^{x}_{x_0} \eta_i 	\right)_{i=1,\ldots b} 
\end{array}
\]
We systematically overline objects that live in the Abelian cover.
We have
\begin{eqnarray*}
\frac{1}{t}	\left(\overline{\Phi}\left(\overline{\gamma (t)}\right)-\overline{\Phi}\left(\overline{\gamma (0)}\right)\right)&=&
\frac{1}{t}\int^{t}_{0}d\overline{\Phi} (\dot{\gamma}(t))dt \\
&=& \frac{1}{t}\left( \int^{t}_{0} \eta_i ( \dot{\gamma}(t)	)\right)_{i=1,\ldots b} \\
&\longrightarrow & \left( \left\langle a_i,h \right\rangle\right)_{i=1,\ldots b}
\end{eqnarray*}
when $t \longrightarrow \pm \infty$, for $\mu$-almost every orbit $\gamma$, by Birkhoff's Ergodic Theorem. In particular the line segment $\overline{\Phi} (\overline{\gamma }(\R))$
is unbounded left and right, thus its projection $\Phi (\gamma (\R))$ to $\T^{b}$ is dense in a subtorus whose  dimension is 
\[
I\left(\left( \left\langle a_i,h \right\rangle\right)_{i=1,\ldots b}\right).
\]

 The latter is precisely  $I(h)$ since $a_1,\ldots a_b$ is an integer basis of $H_1 (M,\R)$.  This proves the claim.

Consequently,  $\overline{\Phi}\left(\overline{\spt \mu}\right)$ is dense in a subspace $D$ of dimension $d$, foliated by straight lines with direction $h$. 
Next we consider
\[
\begin{array}{rcl}
\overline{\Psi} \co \overline{M} & \longrightarrow & \R^{b-1}  \\
x & \longmapsto & \left( \int^{x}_{x_0} \omega'_i 	\right)_{i=1,\ldots b-1} .
\end{array}
\]

\begin{remark}
The reason to lift to the Abelian cover is that the corresponding map from $M$ would take values in a non-separated topological space, thus precluding any discussion of Hausdorff dimension.
\end{remark}
Observe that  the map $\overline{\Psi}$ satisfies a H\H{o}lder condition of order two along $\overline{\Azero}$, and that $\overline{\Psi}= \Lambda \circ \overline{\Phi}$ so $\overline{\Psi}(\overline{\spt \mu})$ is dense in a
subspace of dimension $d-1$. 

Take a Lipschitz vector field $X$ on $M$ such that for every orbit in $\Azero$, for all $t \in \R$, $\dot{\gamma}(t)=X(\gamma(t))$. Such a vector field exists by Mather's Graph Theorem. It does not vanish  on  $\Azero$ because we assume that the zero energy level  does not meet the zero section of $TM$. 

Let $(U_i,f_i),i=1\ldots s$ be such that :
\begin{itemize}
	\item the $U_i$ are open sets of $M$ that cover $\spt \mu$
	\item each $f_i$ is a diffeomorphism from $U_i$ to the open unit ball in $\R^{n}$ 
	\item for all $x\in U_i$, $d_{x}f_i.X(x)=(1,0\ldots ,0)$.
\end{itemize}
In each $U_i$ define
	\[N_i := f^{-1}_{i}\left(\left\{x_1=0\right\}\cap f_i(U_i)\right)
\]
where $x_1$ denotes the first coordinate in $\R^{n}$. That is,  $N_i$ is a transverse section to $X$ in $U_i$. The sets $U_i$ cover $\spt \mu$ so the $\Phi (U_i)$ cover $\Phi (\spt \mu )$, hence one of the $\Phi (U_i\cap \spt \mu)$, say $\Phi (U_1\cap \spt \mu)$, contains an open set of a $d$-dimensional torus. Thus $\overline{\Phi} (\overline{U_1}\cap \overline{\spt \mu})$ contains an open set of  $\R^d$, whence $\overline{\Psi} (\overline{U}_1\cap \overline{\spt \mu})$ contains an open set of  $\R^{d-1}$. But $\overline{\Psi} (\overline{U}_1\cap \overline{\spt \mu})=\overline{\Psi} (\overline{N}_1\cap \overline{\spt \mu})$
because $N_1$ meets every orbit in $U_1$ and $\overline{\Psi}$ is constant on the orbits in $\spt \mu$.

Therefore  the Hausdorff dimension of $\overline{\Psi}\left( \overline{N}_1 \cap \overline{\spt \mu} \right)$ is at least $d-1$, whence by the H\H{o}lder property the Hausdorff dimension of $\overline{N} \cap \overline{\spt \mu}$ is $ \geq 2d-2$. This implies 
$$\dim_H (\overline{\spt \mu)} \geq 2d-1, $$
but recall that $ \overline{\spt \mu}$ is a countable union of lifts of $\spt \mu$ so finally
	\[\dim_H (\spt \mu) \geq 2d-1.
\]
\qed

{\small

\bigskip

\noindent

Math\'ematiques, Universit\'e Montpellier II, France\\
e-mail : massart@math.univ-montp2.fr
}

\end{document}